# A BERNSTEIN–VON MISES THEOREM IN THE NONPARAMETRIC RIGHT-CENSORING MODEL[1]

BY YONGDAI KIM AND JAEYONG LEE

*Seoul National University*

In the recent Bayesian nonparametric literature, many examples have been reported in which Bayesian estimators and posterior distributions do not achieve the optimal convergence rate, indicating that the Bernstein–von Mises theorem does not hold. In this article, we give a positive result in this direction by showing that the Bernstein–von Mises theorem holds in survival models for a large class of prior processes neutral to the right. We also show that, for an arbitrarily given convergence rate $n^{-\alpha}$ with $0 < \alpha \le 1/2$, a prior process neutral to the right can be chosen so that its posterior distribution achieves the convergence rate $n^{-\alpha}$.

**1. Introduction.** The asymptotic properties of posterior distributions and Bayes estimators in nonparametric models have been given much attention in the recent literature. Diaconis and Freedman (1986) opened the discussion in this area by showing that in nonparametric models even an innocent looking prior can produce an inconsistent posterior. This disturbing result stirred Bayesians, because it says that a Bayesian can be more and more sure of a wrong parameter value as the sample size increases. It also initiated research efforts to garner "safe" priors in the asymptotic sense. For the research work regarding posterior consistency, see Freedman (1963), Schwartz (1965), Barron, Schervish and Wasserman (1999) and Ghosal, Ghosh and Ramamoorthi (1999). In the context of survival models, Kim and Lee (2001) showed that not all the prior processes neutral to the right have consistent posterior distributions and gave sufficient conditions for the consistency.

Received February 2001; revised May 2003.

[1]Supported in part by KOSEF through the Statistical Research Center for Complex Systems at Seoul National University.

Cox (1993) and Zhao (2000) showed that this unfortunate phenomenon continues to occur in the posterior convergence rate. For example, Zhao (2000) showed that in an infinite dimensional normal model, there is no independent normal prior supported on the parameter space that has a Bayes estimator that attains the optimal minimax rate. (In the same article, however, she constructed a class of priors, mixtures of normal priors supported on the parameter space, which achieves the optimal minimax rate.) These examples cast doubt on the Bernstein–von Mises theorem in nonparametric models even with the prior that has a consistent posterior.

The Bernstein–von Mises theorem states that the posterior distribution centered at the maximum likelihood estimator (MLE) is asymptotically equivalent to the sampling distribution of the MLE. Due to the recent advent of the Markov chain Monte Carlo method, Bayesians' computational ability exceeds that of frequentists. In the situations where frequentists do not have a computational tool while Bayesians do, frequentists often use the Bayesian credible set as a frequentist confidence interval. The theoretical justification of this practice is the Bernstein–von Mises theorem. Hence, if the Bernstein–von Mises theorem does not hold, this practice is not warranted. The Bernstein–von Mises theorem is squarely important to Bayesians as well, because invalidity of the Bernstein–von Mises theorem often means that a Bayesian credible set has zero efficiency relative to the frequentist confidence interval.

In this article we provide a positive result in this direction by showing that the Bernstein–von Mises theorem does hold in survival models for a large class of prior processes. Indeed, for popular prior processes such as Dirichlet, beta and gamma processes, the Bernstein–von Mises theorem holds. The situation is subtle, however. In an example provided in Section 4, we also show that for any given $0 < \alpha \leq 1/2$, there is a consistent prior process neutral to the right that has a posterior convergence rate that is exactly $n^{-\alpha}$. This result suggests that, for a given model and data, one prior process can be much slower extracting information from the data than another. This confirms the findings in the literature that posterior consistency does not guarantee the optimal convergence rate and in practice a prior must be carefully examined before it is used. In the same example, an interesting prior process is found. This prior process achieves the optimal posterior convergence rate, but its posterior distribution is not equivalent to the sampling distribution of the MLE; hence, the Bernstein–von Mises theorem does not hold. This example shows that the optimal convergence rate does not guarantee the Bernstein–von Mises theorem.

The Bernstein–von Mises theorem for parametric models is a well-known result. See, for instance, Section 7.4.2 of Schervish (1995) and references therein. Previous research on the Bernstein–von Mises theorem for nonparametric models includes Lo (1983, 1986, 1993), Brunner and Lo (1996),

Diaconis and Freedman ([1998](#)), Conti ([1999](#)) and Freedman ([1999](#)). Among them, Lo ([1983](#), [1986](#), [1993](#)), Brunner and Lo ([1996](#)) and Conti ([1999](#)) reported some of the earlier positive results on the Bernstein–von Mises theorem for some nonparametric models. See also Ghosal, Ghosh and van der Vaart ([2000](#)) and Shen and Wasserman ([2001](#)) for a related theory of posterior convergence rates.

In Section 2 the survival model and prior processes neutral to the right are briefly introduced. In Section 3 the main result of this article, the Bernstein–von Mises theorem of survival models, is given. In Section 4 a class of prior processes with arbitrary posterior convergence rate $n^{-\alpha}$, $0 < \alpha \le 1/2$, and a simulation study are given. The proof of the Bernstein–von Mises theorem is given in Section 5.

**2. Survival models and processes neutral to the right.** Let $X_1, \ldots, X_n$ be i.i.d. survival times with cumulative distribution function (c.d.f.) $F$ and let $C_1, \ldots, C_n$ be independent censoring times with c.d.f. $G$, independent of the $X_i$'s. Since the observations are subject to right censoring, we observe only $(T_1, \delta_1), \ldots, (T_n, \delta_n)$, where $T_i = \min(C_i, X_i)$ and $\delta_i = I(X_i \le C_i)$. Let $D_n = \{(T_1, \delta_1), \ldots, (T_n, \delta_n)\}$. Let $A$ be the cumulative hazard function (c.h.f.) of $F$, $A(t) = \int_0^t dF(s)/(1 - F(s-))$.

We say that a prior process on c.d.f. $F$ is a process neutral to the right if the corresponding c.h.f. $A$ is a nonstationary subordinator (a positive nondecreasing independent increment process) such that $A(0) = 0$, $0 \le \Delta A(t) \le 1$ for all $t$ with probability 1 and either $\Delta A(t) = 1$ for some $t > 0$ or $\lim_{t \to \infty} A(t) = \infty$ with probability 1. See Doksum ([1974](#)) for the original definition of processes neutral to the right and see Hjort ([1990](#)), Kim ([1999](#)) and Kim and Lee ([2001](#)) for the connection between the definition given here and Doksum's definition. In what follows, the term *subordinator* is used for a prior process of c.h.f. $A$ which induces a process neutral to the right on $F$.

Kim ([1999](#)) used the following characterization of subordinators. This characterization can be dated back to Lévy [see the note in Breiman ([1968](#)), page 318]. Similar characterization can also be found in Theorem 6.3VIII in Daley and Vere-Jones ([1988](#)) and Theorem 3 in Fristedt and Gray [([1997](#)), page 606]. For any given subordinator $A(t)$ on $[0, \infty)$, there exists a unique random measure $\mu$ on $[0, \infty) \times [0, 1]$ such that

$$(1) \qquad A(t) = \int_{[0,t] \times [0,1]} x \mu(ds, dx).$$

In fact, $\mu$ is defined by

$$\mu([0, t] \times B) = \sum_{s \le t} I(\Delta A(s) \in B)$$

for any Borel subset $B$ of $[0, 1]$ and for all $t > 0$. Since $\mu$ is a Poisson random measure [Jacod and Shiryaev (1987), page 70], there exists a unique $\sigma$-finite measure $\nu$ on $[0, \infty) \times [0, 1]$ such that

$$(2) \qquad \qquad \mathrm{E}(\mu([0, t] \times B)) = \nu([0, t] \times B)$$

for all $t > 0$. Conversely, for a given $\sigma$-finite measure $\nu$ such that

$$\int_0^t \int_0^1 x \nu(ds, dx) < \infty$$

for all $t$, there exists a unique Poisson random measure $\mu$ on $[0, \infty) \times [0, 1]$ which satisfies (2) [Jacod (1979)] and so we can construct a subordinator $A$ through (1). Conclusively, we can use $\nu$ to characterize a subordinator $A$.

Suppose that a given subordinator $A$ has fixed discontinuity points at $t_1 < t_2 < \cdots$ and that the Lévy formula is given by

$$\mathrm{E}(\exp(-\theta A(t))) = \left[ \prod_{t_i \le t} \mathrm{E}(\exp(-\theta \Delta A(t_i))) \right] \exp\left( - \int_0^1 (1 - e^{-\theta x}) \, dL_t(x) \right),$$

where $L_t(x)$ is the Lévy measure. Then it can be shown [see Theorem II.4.8 in Jacod and Shiryaev (1987)] that

$$\nu([0, t] \times B) = \int_B dL_t(x) + \sum_{t_i \le t} \int_B dH_i(x)$$

for all $t > 0$ and for any Borel set $B$ of $[0, 1]$, where $H_i(x)$ is the distribution function of $\Delta A(t_i)$. When there are no fixed discontinuities, $\mu$ is a Poisson random measure defined on $[0, \infty) \times [0, 1]$ with intensity measure $\nu$ and $dL_t(x) = \int_{[0,t]} \nu(ds, dx)$. Hence, the measure $\nu$ simply extends $dL_t$ by incorporating the fixed discontinuity points. However, this simple extension provides a convenient notational device. The posterior distribution, which typically has many fixed discontinuity points, can be summarized neatly by use of the corresponding measure $\nu$ without separating out the stochastically continuous part and the fixed discontinuity points as was done in previous work [Ferguson and Phadia (1979) and Hjort (1990)]. For this reason, we call $\nu$ simply the *Lévy measure* of $A$.

From the Lévy measure $\nu$, we can easily calculate the mean and variance of the subordinator using the formulas [Kim (1999)]

$$(3) \qquad \qquad \mathrm{E}(A(t)) = \int_0^t \int_0^1 x \nu(ds, dx)$$

and

$$\mathrm{Var}(A(t)) = \int_0^t \int_0^1 x^2 \nu(ds, dx) - \sum_{s \le t} \left( \int_0^1 x \nu(\{s\}, dx) \right)^2.$$

These formulas constitute basic facts for the asymptotic theory of the posterior and will be used subsequently herein.

The characterization of subordinators with Lévy measures is also convenient in representing the posterior distribution, for the class of processes neutral to the right is conjugate with respect to right censored survival data. Suppose a priori $A$ is a subordinator with Lévy measure

$$\nu(ds, dx) = f_s(x) \, dx \, ds \qquad \text{for } s \geq 0 \text{ and } 0 \leq x \leq 1, \tag{4}$$

with $\lim_{t \to \infty} \int_0^t \int_0^1 x f_s(x) \, dx \, ds = \infty$. Then the posterior distribution of $A$ given $D_n$ is again a subordinator with Lévy measure $\nu^p$ given by

$$\nu^p(ds, dx) = (1-x)^{Y_n(s)} f_s(x) \, dx \, ds + dH_s(x) \frac{1}{\Delta N_n(s)} \, dN_n(s), \tag{5}$$

where $H_s(x)$ is a distribution function on $[0, 1]$ and is defined by

$$dH_s(x) \propto x^{\Delta N_n(s)} (1-x)^{Y_n(s) - \Delta N_n(s)} f_s(x) \, dx$$

and $N_n(t) = \sum_{i=1}^n I(T_i \leq t, \delta_i = 1)$, $Y_n(t) = \sum_{i=1}^n I(T_i \geq t)$, $\Delta N_n(t) = N_n(t) - N_n(t-)$. Note that the posterior process is the sum of stochastically continuous and discrete parts, which correspond to the first and the second terms in (5), respectively. Note also that $H_s$ is the distribution of jump size at $s$ if $\Delta N_n(s) \neq 0$. This fact is used later. For the proof of (5), see Hjort (1990) or Kim (1999).

Let $F_0$ be the true distribution of the $X_i$'s and let $A_0$ be the c.h.f. of $F_0$. We will study the asymptotic behavior of $A$ on a fixed compact interval $[0, \tau]$. Throughout this article we assume the following two conditions:

CONDITION C1. $F_0(\tau-) < 1$ and $G(\tau-) < 1$.

CONDITION C2. $A_0$ is continuous on $[0, \tau]$.

Condition C1 guarantees that $Y_n(\tau) \to \infty$ as $n \to \infty$ with probability 1, which is essential for the asymptotic theory of survival models. Condition C2 implies that $\Delta N_n(s)$ has a value of either 0 or 1.

**3. Bernstein–von Mises theorem.** Assume that a priori $A$ is a nonstationary subordinator with Lévy measure

$$\nu([0, t] \times B) = \int_0^t \int_B \frac{1}{x} g_s(x) \, dx \, \lambda(s) \, ds, \tag{6}$$

where $\int_0^1 g_t(x) \, dx = 1$ for all $t \in [0, \tau]$.

REMARK. Comparing (4) and (6), we can see that $\lambda(t) = \int_0^1 x f_t(x) \, dx$ and $g_t(x) = x f_t(x) / \lambda(t)$ provided $\lambda(t) > 0$.

We need the following conditions for the Bernstein–von Mises theorem:

CONDITION A1.    $g^* = \sup_{t \in [0,\tau], x \in [0,1]} (1-x) g_t(x) < \infty$.

CONDITION A2.    There exists a function $q(t)$ defined on $[0, \tau]$ such that $0 < \inf_{t \in [0,\tau]} q(t) \le \sup_{t \in [0,\tau]} q(t) < \infty$ and, for some $\alpha > 0$ and $\varepsilon > 0$,

$$\sup_{t \in [0,\tau], x \in [0,\varepsilon]} \left| \frac{g_t(x) - q(t)}{x^\alpha} \right| < \infty.$$

CONDITION A3.    $\lambda(t)$ is bounded and positive on $(0, \tau)$.

The convergence rate of the posterior distribution depends mainly on the behavior of the prior process in the neighborhood of 0. This is because the jump sizes of the posterior process get smaller as $n$ gets larger. Condition A1 is a technical one to make the posterior mass of the jump sizes of the fixed discontinuity points outside the neighborhood of 0 be asymptotically negligible. Condition A2 is the main condition, in which $\alpha$ measures the smoothness of $g_t(x)$ in $x$ around 0. The constant $\alpha$ plays a crucial role in determining the convergence rate of the posterior distribution. In fact, the Bernstein–von Mises theorem may not hold if $\alpha \le 1/2$. For an example, see Section 4. The boundedness of $\lambda$ in Condition A3 makes the posterior distribution eventually be dominated by data. The positiveness of $\lambda$ in Condition A3 is also necessary. Suppose $\lambda(t) = 0$ for $t \in [c, d]$, where $0 < c < d < \tau$. Then both the prior and posterior put mass 1 to the set of c.h.f.s $A$, with $A(d) = A(c)$. Hence the posterior distribution of $A(d) - A(c)$ has mass 1 at 0 and the Bernstein–von Mises theorem does not hold unless $A_0(d) = A_0(c)$.

Before stating the theorems, we introduce some notation. For a given random variable $Z_n$, we write $Z_n = O(n^\delta)$ with probability 1 if there exists a constant $M > 0$ such that $|Z_n|/n^\delta \le M$ for all but finitely many $n$ with probability 1. Let $\delta_a$ be the degenerate probability measure at $a$. Denote by $\mathcal{L}(X|Y)$ the conditional distribution of $X$ given $Y$. Let $W$ be a standard Brownian motion and let $\hat{A}_n$ be the Aalen–Nelson estimator defined by $\hat{A}_n(t) = \int_0^t dN_n(s)/Y_n(s)$. The sampling distribution of $\sqrt{n}(\hat{A}_n - A_0)$ converges in distribution to $W(U_0(\cdot))$, where $U_0(t) = \int_0^t dA_0(s)/Q(s)$, with $Q(t) = \Pr(T_1 \ge t)$ [see Theorem IV.1.2 in Andersen, Borgan, Gill and Keiding (1993)]. Here $U_0$ is well defined, because $\inf_{t \in [0,\tau]} Q(t) > 0$ due to Condition C1.

The following theorem is a general result on the convergence of the posterior distribution. The Bernstein–von Mises theorem and an example of suboptimal convergent rates in Section 4 will be based on this theorem. Let $q_n$ be the number of distinct uncensored observations and let $t_1 < t_2 < \cdots < t_{q_n}$ be the distinct uncensored observations. Let $A_d(t) = \sum_{i=1}^{q_n} \Delta A(t_i)$. Let $D[0, \tau]$ be the space of cadlag functions on $[0, \tau]$ equipped with the uniform topology and the ball $\sigma$-field.

Theorem 1. *Under Conditions* A1–A3*:*

(i) $\mathcal{L}(\sqrt{n}(A(\cdot) - A_d(\cdot))|\mathrm{D}_n) \xrightarrow{d} \delta_0$ *on* $D[0, \tau]$ *with probability* 1*;*

(ii) $\mathcal{L}(\sqrt{n}(A_d(\cdot) - \mathrm{E}(A_d(\cdot)|\mathrm{D}_n))|\mathrm{D}_n) \xrightarrow{d} W(U_0(\cdot))$ *on* $D[0, \tau]$ *with probability* 1*;*

(iii) $\sup_{t \in [0, \tau]} |\mathrm{E}(A_d(t)|\mathrm{D}_n) - \hat{A}_n(t)| = O(n^{-\min\{1, \alpha\}})$ *with probability* 1.

The proof is given in Section 5.

Part (i) of Theorem 1 states that the stochastically continuous part of the posterior process, $A - A_d$, vanishes with a rate faster than the optimal rate, $n^{-1/2}$. Part (ii) states that the fixed discontinuous part of the posterior process, $A_d$, centered at its mean is asymptotically equivalent to the frequentist sampling distribution of $\hat{A}_n$ since $W(U_0(t))$ in Theorem 1(ii) is the limiting sampling distribution of $\sqrt{n}(\hat{A}_n(t) - A_0(t))$. Part (iii) states that the difference of the posterior mean of $A_d$ and $\hat{A}_n$ vanishes with varying order, $n^{-\min\{1, \alpha\}}$, for $\alpha > 0$. Hence, if $\alpha < 1/2$, the overall convergence rate of the posterior distribution could be dominated by the convergence rate of (iii), which results in suboptimal convergence rates. Indeed, in Section 4 such an example is given.

Although a rigorous proof of Theorem 1 is given in Section 5, we sketch the proof here. For (i), we first approximate the first two moments of the posterior distribution of $A$ by those of the posterior with a beta process prior (see Example 1 for a definition of beta process). Since the closed forms of the first two moments of the posterior with the beta process prior are known [Hjort (1990)], one can easily prove (i) using Lemma 7. Part (iii) is proved similarly. For (ii), the posterior distribution of $A_d$ consists of the sum of independent random variables, and so the central limit theorem for independent random variables [e.g., Theorem 19 in Section V.4 in Pollard (1984)] can be applied.

Theorem 2 (Bernstein–von Mises theorem). *Under Conditions* A1–A3 *with* $\alpha > 1/2$,

$$\mathcal{L}(\sqrt{n}(A(\cdot) - \hat{A}_n(\cdot))|\mathrm{D}_n) \xrightarrow{d} W(U_0(\cdot))$$

*on* $D[0, \tau]$ *with probability* 1.

Proof. This theorem is an immediate consequence of Theorem 1, because we can decompose

$$n^{1/2}(A(t) - \hat{A}(t)) = n^{1/2}(A(t) - A_d(t)) + n^{1/2}(A_d(t) - \mathrm{E}(A_d(t)|\mathrm{D}_n))$$
$$+ n^{1/2}(\mathrm{E}(A_d(t)|\mathrm{D}_n) - \hat{A}_n(t)). \qquad \square$$

COROLLARY 1. *Under the same conditions in Theorem* 2,

$$\mathcal{L}(\sqrt{n}(S(\cdot) - \hat{S}_n(\cdot))|\mathrm{D}_n) \xrightarrow{d} -S_0(\cdot)W(U_0(\cdot))$$

*on* $D[0, \tau]$ *with probability* 1, *where* $S, \hat{S}_n$ *and* $S_0$ *are the corresponding survival functions of* $A, \hat{A}_n$ *and* $A_0$.

PROOF. Note that the survival function is recovered from the cumulative hazard function by the product integration operator which is Hadamard differentiable. The result follows from the functional delta method. See Gill (1989).  □

REMARK. If Conditions A1–A3 as well as Condition C1 hold for all $\tau > 0$, Theorem 2 and Corollary 1 are valid on $D[0, \infty)$, because the weak convergence on $D[0, \infty)$ is defined by the weak convergence on $D[0, \tau]$ for all $\tau > 0$ [Pollard (1984)].

A convenient sufficient condition for Condition A2 with $\alpha = 1$ can be given as follows. Suppose that for some $\varepsilon > 0$

$$\sup_{t \in [0, \tau], x \in (0, \varepsilon)} |g_t^{(1)}(x)| < \infty, \tag{7}$$

where $g_t^{(1)}(x)$ is the first derivative of $g_t(x)$ in $x$ on $[0, 1]$. Then, by the mean value theorem, Condition A2 holds with $\alpha = 1$ and $q(t) = g_t(0)$.

In the next three examples, we illustrate that the Bernstein–von Mises theorem holds for beta, Dirichlet and gamma prior processes.

EXAMPLE 1 (Beta processes). The beta process with mean $\Lambda$ and scale parameter $c$ is a nonstationary subordinator with Lévy measure $\nu$, $\nu(dt, dx) = c(t)x^{-1}(1-x)^{c(t)-1} \, dx \, d\Lambda(t)$. Suppose $\Lambda(t) = \int_0^t \lambda(s) \, ds$, where $\lambda(t)$ is positive continuous on $(0, \tau)$ and $0 < \inf_{t \in [0, \tau]} c(t)(= c_*) \le \sup_{t \in [0, \tau]} c(t)(= c^*) < \infty$. Condition A1 is true because

$$\sup_{t \in [0, \tau], x \in [0, 1]} |(1-x)g_t(x)| = \sup_{t \in [0, \tau], x \in [0, 1]} |c(t)(1-x)^{c(t)}| \le c^* < \infty.$$

For Condition A2, since $g_t^{(1)}(x) = c(t)(c(t) - 1)(1-x)^{c(t)-2}$, we have

$$\sup_{t \in [0, \tau], x \in (0, \varepsilon)} |g_t^{(1)}(x)| \le c^*(c^* + 1) \max\{1, (1-\varepsilon)^{c_*-2}\}.$$

Thus, by (7), Condition A2 holds with $q(t) = c(t)$. Since Condition A3 is assumed, the Bernstein–von Mises theorem holds.

EXAMPLE 2 (Dirichlet processes). Hjort (1990) showed that when the prior of the distribution $F$ is the Dirichlet process with base measure $\alpha$, the induced prior of the c.h.f. is the beta process with $c(t) = \alpha([0, \infty))(1 - H(t))$ and $\Lambda(t)$, the c.h.f. of $H(t)$, where $H(t) = \alpha([0, t])/\alpha([0, \infty))$. Suppose $\Lambda(t) = \int_0^t \lambda(s)\, ds$. Then if $\lambda(t)$ is positive bounded on $(0, \tau)$ and $H(\tau) < 1$, then, as in Example 1, it can be shown that Conditions A1–A3 are satisfied.

EXAMPLE 3 (Gamma processes). A priori, assume that $Y(t) = -\log(1 - F(t))$ is a gamma process with parameters $(H(t), d(t))$ with $H(t) = \int_0^t h(s)\, dx$. Here the gamma process with parameters $(H(t), d(t))$ is defined by $Y(t) = \int_0^t \frac{1}{d(s)}\, dX(s)$, where $X(t)$ is a subordinator that has a marginal distribution of $X(t)$ that is a gamma distribution with parameters $(\int_0^t d(s)\, dH(s), 1)$. See Lo (1982) for details. This prior process was used by Doksum (1974), Kalbfleisch (1978) and Ferguson and Phadia (1979) . Since

$$\log \mathrm{E}(\exp(-\theta Y(t))) = \int_0^t \int_0^\infty (e^{-\theta x} - 1)\frac{d(s)}{x}\exp(-d(s)x)\, dx\, dH(s),$$

it can be shown that the c.h.f. $A$ of $F$ is a subordinator with Lévy measure $\nu$ given by

$$\nu([0, t] \times B) = \int_0^t c(s)\int_B \frac{1}{-\log(1 - x)}(1 - x)^{d(s) - 1}\, dx\, d\Lambda(s),$$

where

$$c(t) = \left(\int_0^1 \frac{x}{-\log(1 - x)}(1 - x)^{d(t) - 1}\, dx\right)^{-1}$$

and

$$\Lambda(t) = \int_0^t \frac{d(s)}{c(s)}\, dH(s).$$

Therefore, we have

$$g_t(x) = c(t)\frac{x}{-\log(1 - x)}(1 - x)^{d(t) - 1}, \qquad 0 \le x \le 1,$$

and $\lambda(t) = d(t)h(t)/c(t)$.

Suppose $h(t)$ is positive and bounded on $t \in (0, \tau)$ and $0 < \inf_{t \in [0, \tau]} d(t)(= d_*) \le \sup_{t \in [0, \tau]} d(t)(= d^*) < \infty$. We will show that Conditions A1–A3 hold under these conditions. First, we show that $0 < \inf_{t \in [0, \tau]} c(t)(= c_*) \le \sup_{t \in [0, \tau]} c(t)(= c^*) < \infty$. Note that

$$\inf_{t \in [0, \tau]} c(t) = \left(\int_0^1 \frac{x}{-\log(1 - x)}(1 - x)^{d_* - 1}\, dx\right)^{-1}$$

$$= \left(\int_0^1 \frac{x(1 - x)^{d_*/2}}{-\log(1 - x)}(1 - x)^{d_*/2 - 1}\, dx\right)^{-1} \ge \left(\frac{m}{d_*/2}\right)^{-1} > 0,$$

where $m = \sup_{x \in [0,1]} -x(1-x)^{d_*/2}/\log(1-x)$. By a similar argument, we can show that $\sup_{t \in [0,\tau]} c(t) < \infty$. Now, Condition A1 follows because

$$\sup_{t \in [0,\tau], x \in [0,1]} |(1-x)g_t(x)| = \sup_{t \in [0,\tau], x \in [0,1]} \left| c(t)\frac{x}{-\log(1-x)}(1-x)^{d(t)} \right|$$

$$\leq \sup_{t \in [0,\tau], x \in [0,1]} \left| c(t)\frac{x}{-\log(1-x)}(1-x)^{d_*} \right|$$

$$\leq c^* \sup_{x \in [0,1]} \left| \frac{x(1-x)^{d_*}}{-\log(1-x)} \right| < \infty.$$

Similarly, Condition A2 can be shown by (7) with $q(t) = c(t)$ and Condition A3 follows from $\inf_{t \in [0,\tau]} c(t) > 0$.

**4. An example: suboptimal convergence rates.** In this section, we show that, for a given convergence rate $n^{-\alpha}$ with $0 < \alpha \leq 1/2$, there exists a prior process neutral to the right whose posterior convergence rate is $n^{-\alpha}$. Consider the class of prior processes neutral to the right with Lévy measure

$$(8) \qquad \nu_\alpha(dt, dx) = \frac{1}{x}(1 + x^\alpha)\, dx\, dt, \qquad x \in (0,1], t \geq 0.$$

In the next theorem we show that, for each $0 < \alpha \leq 1/2$, the posterior with the prior process $\nu_\alpha$ achieves convergence rate $n^{-\alpha}$.

THEOREM 3. *A priori let $A$ be a subordinator with Lévy measure $\nu_\alpha$ in* (8). *Then:*

(i) *For $0 < \alpha < 1/2$, $\mathcal{L}(n^\alpha(A(\cdot) - \hat{A}_n(\cdot))|\mathrm{D}_n) \xrightarrow{d} \delta_{J_\alpha(\cdot)}$ on $D[0,\tau]$ with probability 1, where $J_\alpha(t) = \alpha\Gamma(\alpha+1)\int_0^t dA_0(s)/Q^\alpha(s)$.*

(ii) *For $\alpha = 1/2$, $\mathcal{L}(n^{1/2}(A(\cdot) - \hat{A}_n(\cdot))|\mathrm{D}_n) \xrightarrow{d} W(U_0(\cdot)) + J_{1/2}(\cdot)$ on $D[0,\tau]$ with probability 1.*

(iii) *For $\alpha > 1/2$, $\mathcal{L}(n^{1/2}(A(\cdot) - \hat{A}_n(\cdot))|\mathrm{D}_n) \xrightarrow{d} W(U_0(\cdot))$ on $D[0,\tau]$ with probability 1.*

REMARK 1. When $0 < \alpha < 1/2$, the posterior convergence rate is $n^{-\alpha}$, which is slower than the optimal rate $n^{-1/2}$.

REMARK 2. When $\alpha = 1/2$, the posterior convergence rate is optimal, but the limiting posterior distribution is the limiting sampling distribution of the Aalen–Nelson estimator plus a bias term $J_{1/2}$. So the Bayesian credible set does not have appropriate frequentist coverage probability, although it has the optimal posterior convergence rate.

REMARK 3. The Bernstein–von Mises theorem holds when $\alpha > 1/2$. Although we do not know whether Conditions A1–A3 are necessary and sufficient conditions for the Bernstein–von Mises theorem, this example shows that these conditions are fairly minimal.

To prove Theorem 3 we need the following lemma, the proof of which can be found in Appendix A.

LEMMA 1. *For* $0 < \alpha \leq 1/2$,

$$\sup_{t \in [0,\tau]} |n^\alpha (\mathrm{E}(A_d(t)|\mathrm{D}_n) - \hat{A}_n(t)) - J_\alpha(t)| \to 0$$

*with probability* 1.

PROOF OF THEOREM 3. It is easy to see that $\nu_\alpha$ in (8) satisfies Conditions A1–A3 with $q(t) = (\alpha + 1)/(\alpha + 2)$ and $\lambda(t) = 1$. Now note that

$$n^\alpha (A(t) - \hat{A}_n(t)) = n^\alpha (A(t) - A_d(t)) + n^\alpha (A_d(t) - \mathrm{E}(A_d(t)|\mathrm{D}_n))$$
$$+ n^\alpha (\mathrm{E}(A_d(t)|\mathrm{D}_n) - \hat{A}_n(t)).$$

The first term of the right-hand side converges weakly to 0 for all $\alpha > 0$; the second term converges weakly to 0 for $\alpha < 1/2$ and converges weakly to $W(U_0(\cdot))$ for $\alpha \geq 1/2$ by Theorem 1. Finally, the third term converges weakly to $J_\alpha$ for $\alpha \leq 1/2$ and converges weakly to 0 for $\alpha > 1/2$ by Lemma 1. Hence, the proof is complete by Slutsky's theorem. □

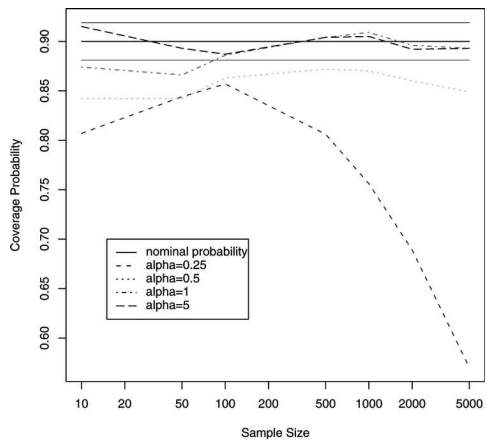

FIG. 1. *Empirical coverage probabilities of the Bayesian credible set of* $A(2)$, *the c.h.f. at* $t = 2$, *with nominal level* 90%. *Empirical coverage probabilities are based on* 1000 *data sets for each of sample sizes* $n = 10, 50, 100, 500, 1000, 2000, 5000$ *with the prior* (8) *at* $\alpha = 0.25, 0.5, 1, 5$. *The three solid lines represent the nominal level and* 2 *standard errors away from it. The dotted lines are the empirical coverage probabilities.*

Theorem 3 shows that the posterior with the prior (8) with $\alpha > 1/2$ can be used to construct an asymptotically valid frequentist confidence interval, while the posterior with $\alpha \leq 1/2$ cannot. A simulation study was conducted to see the effect of $\alpha$ and sample size $n$ on empirical coverage probability. Right censored data were generated from Exponential(1) for the survival time and Exponential(0.25) for the right-censoring times, which amounts to censoring probability 0.2. For each of seven sample sizes $n = 10, 50, 100, 500, 1000, 2000, 5000,$ 1000 data sets were generated. The posterior distribution was computed, based on an algorithm modified from Lee and Kim (2004), for each data set with the prior (8) for $\alpha = 0.25, 0.5, 1, 5$. The empirical coverage probability is the proportion of the data sets that have credible sets of $A(2)$, the c.h.f. at $t = 2$, that contain the true value $A(2) = 2$. The simulation result is reported in Figure 1. The three solid lines represent the nominal coverage probability 0.9 and 2 standard errors $2\sqrt{0.9 \cdot 0.1/1000} = 0.01897$ away from it. The coverage probability with $\alpha = 0.25$ gets worse as the sample size grows. When $\alpha = 0.5$, the coverage probability shows a difference from the nominal level which does not get smaller as the sample size increases. However, with $\alpha = 1$ and 5, the coverage probability is inside the error bounds from $n = 100$ on. All of these agree with Theorem 3.

**5. Proof of Theorem 1.** Throughout this section, the statements of Theorem 1 are assumed. Let $B(a, b) = \int_0^1 x^{a-1}(1-x)^{b-1}\,dx$. Then Stirling's formula yields that, for $\alpha > 0$,

$$(9) \qquad \lim_{n \to \infty} n^\alpha B(\alpha, n) = \Gamma(\alpha).$$

LEMMA 2. *Let $W_n$ be a sequence of nonnegative stochastic processes on $[0, \tau]$ such that*

$$(10) \qquad \sup_{t \in [0, \tau]} |W_n(t)/n - Q(t)| \to 0$$

*with probability* 1. *Then*

$$\sup_{t \in [0, \tau]} |n^k B(k, W_n(t)) - \Gamma(k)/Q^k(t)| \to 0$$

*with probability* 1 *as $n \to \infty$, for every integer $k \geq 1$.*

PROOF. We can write

$$n^k B(k, W_n(t)) = \left(\frac{n}{W_n(t)}\right)^k W_n^k(t) B(k, W_n(t)).$$

Since (10) implies $\inf_{t\in[0,\tau]} W_n(t) \to \infty$ with probability 1, (9) yields

$$\sup_{t\in[0,\tau]} |W_n^k(t)B(k, W_n(t)) - \Gamma(k)| \to 0$$

with probability 1. Also (10) implies

$$\sup_{t\in[0,\tau]} \left|\left(\frac{n}{W_n(t)}\right)^k - \frac{1}{Q^k(t)}\right| \to 0$$

with probability 1, which completes the proof. $\square$

Let $Y_n^+(t) = Y_n(t) - \Delta N_n(t)$ and $C_k(t) = \int_0^1 x^k(1-x)^{Y_n^+(t)} g_t(x)\,dx$.

LEMMA 3. *For $k \geq 0$,*

$$\sup_{t\in[0,\tau]} |C_k(t) - q(t)B(k+1, Y_n^+(t)+1)| = O(n^{-(k+1+\alpha)}) \tag{11}$$

*and*

$$\sup_{t\in[0,\tau]} |n^{k+1}C_k(t) - q(t)\Gamma(k+1)Q^{-(k+1)}(t)| \to 0 \tag{12}$$

*with probability 1.*

PROOF. For (11), let $p_t(x) = (1-x)(g_t(x) - q(t))/x^\alpha$. Then Conditions A1 and A2 together imply $\sup_{t\in[0,\tau],x\in[0,1]} |p_t(x)|(=p^*) < \infty$. Now

$$|C_k(t) - q(t)B(k+1, Y_n^+(t)+1)|$$

$$= \left|\int_0^1 x^{k+\alpha}(1-x)^{Y_n^+(t)-1} p_t(x)\,dx\right|$$

$$\leq p^* B(k+\alpha+1, Y_n^+(t)).$$

Since $\sup_{t\in[0,\tau]} |Y_n^+(t)/n - Q(t)| \to 0$ with probability 1, Lemma 2 yields

$$\sup_{t\in[0,\tau]} |C_k(t) - q(t)B(k+1, Y_n^+(t)+1)| = O(n^{-(k+1+\alpha)})$$

with probability 1. Equation (12) is an easy consequence of (11) and Lemma 2. $\square$

LEMMA 4. *We have*

$$\sup_{t\in[0,\tau]} \left|\frac{C_k(t)}{C_0(t)} - \frac{B(k+1, Y_n^+(t)+1)}{B(1, Y_n^+(t)+1)}\right| = O(n^{-(k+\alpha)}) \tag{13}$$

*and*

$$\sup_{i=1,\dots,q_n} \mathrm{E}((\Delta A_d(t_i))^k | \mathrm{D}_n) = O(n^{-k}). \tag{14}$$

PROOF.   For (13), we can write

$$\left| \frac{C_k(t)}{C_0(t)} - \frac{B(k+1, Y_n^+(t)+1)}{B(1, Y_n^+(t)+1)} \right|$$

$$(15) \qquad \leq \left| \frac{C_k(t) - q(t)B(k+1, Y_n^+(t)+1)}{C_0(t)} \right|$$

$$(16) \qquad + \left| \frac{B(k+1, Y_n^+(t)+1)(C_0(t) - q(t)B(1, Y_n^+(t)+1))}{C_0(t)B(1, Y_n^+(t)+1)} \right|.$$

Since (12) yields

$$(17) \qquad \inf_{t\in[0,\tau]} nC_0(t) \to \inf_{t\in[0,\tau]} q(t)Q^{-1}(t) > 0,$$

(11) implies

$$n^{k+\alpha} \sup_{t\in[0,\tau]} \left| \frac{C_k(t) - q(t)B(k+1, Y_n^+(t)+1)}{C_0(t)} \right|$$

$$\leq \frac{n^{k+\alpha+1} \sup_{t\in[0,\tau]} |C_k(t) - q(t)B(k+1, Y_n^+(t)+1)|}{\inf_{t\in[0,\tau]} nC_0(t)}$$

$$= O(1)$$

with probability 1, and hence (15) is $O(n^{-(k+\alpha)})$. On the other hand, since (11) yields

$$\inf_{t\in[0,\tau]} nB(1, Y_n^+(t)+1) \to \inf_{t\in[0,\tau]} Q^{-1}(t) > 0,$$

(11) together with Lemma 2 and (17) implies

$$n^{k+\alpha} \sup_{t\in[0,\tau]} \left| \frac{B(k+1, Y_n^+(t)+1)(C_0(t) - q(t)B(1, Y_n^+(t)+1))}{C_0(t)B(1, Y_n^+(t)+1)} \right|$$

$$\leq \frac{\sup_{t\in[0,\tau]} |n^{k+1}B(k+1, Y_n^+(t)+1)|}{\inf_{t\in[0,\tau]} nC_0(t)}$$

$$\times \frac{\sup_{t\in[0,\tau]} |n^{1+\alpha}(C_0(t) - q(t)B(1, Y_n^+(t)+1))|}{\inf_{t\in[0,\tau]} nB(1, Y_n^+(t)+1)}$$

$$= O(1)$$

with probability 1, and so (16) is $O(n^{-(k+\alpha)})$, which completes the proof of (13).

For (14), note that the distribution function of $\Delta A_d(t_i)$ is $H_{t_i}(x)$. Hence $\mathrm{E}((\Delta A_d(t_i))^k | \mathrm{D}_n) = C_k(t_i)/C_0(t_i)$, which together with (14) completes the proof.   □

Proof of Theorem 1(i). Since a posteriori $A - A_d$ is a Lévy process with Lévy measure $\nu_c$ given by $\nu_c(dt, dx) = x^{-1}(1-x)^{Y_n(t)} g_t(x) \, dx \, \lambda(t) \, dt$, Condition A1 and Lemma 2 with (3) imply

$$\mathrm{E}(A(t) - A_d(t)|\mathrm{D}_n) \leq g^* \int_0^\tau \lambda(s) \, ds B(1, Y_n(\tau)) = O(n^{-1})$$

with probability 1. Similarly,

$$V(A(t) - A_d(t)|\mathrm{D}_n) \leq g^* \int_0^\tau \lambda(s) \, ds B(2, Y_n(\tau)) = O(n^{-2})$$

with probability 1. Hence the proof is completed by Lemma 7 (in Appendix B). □

Proof of Theorem 1(ii). Let $Z_n(t) = \sqrt{n}(A_d(t) - \mathrm{E}(A_d(t)|\mathrm{D}_n))$. Since $Z_n$ is a Lévy process, we can utilize Theorem 19 in Section V.4 in Pollard (1984). We first prove the convergence of finite dimensional distributions by showing Lyapounov's condition. Suppose $0 \leq s < t \leq \tau$ are given. Note that

$$Z_n(t) - Z_n(s) = \sum_{s < t_i \leq t} \sqrt{n}(\Delta A_d(t_i) - \mathrm{E}(\Delta A_d(t_i)|\mathrm{D}_n)).$$

Let

$$\sup_{i=1,\dots,q_n} \mathrm{E}\Big[\Big(\sqrt{n}(\Delta A_d(t_i) - \mathrm{E}(\Delta A_d(t_i)|\mathrm{D}_n))\Big)^4 |\mathrm{D}_n\Big] = V_n.$$

Then (14) in Lemma 4 implies $V_n = O(n^{-2})$ with probability 1. Because $\sup_{t \in [0,\tau]} \int_0^t dN_n(u) = O(n)$,

$$(18) \quad \begin{aligned} \sum_{s < t_i \leq t} \mathrm{E}\Big[\Big(\sqrt{n}(\Delta A_d(t_i) - \mathrm{E}(\Delta A_d(t_i)|\mathrm{D}_n))\Big)^4 |\mathrm{D}_n\Big] \\ \leq \int_s^t V_n \, dN_n(u) \to 0 \end{aligned}$$

with probability 1.

On the other hand, let

$$\begin{aligned} W_{ni} = \mathrm{E}((\Delta A_d(t_i))^2|\mathrm{D}_n) - (\mathrm{E}(\Delta A_d(t_i)|\mathrm{D}_n))^2 \\ - \left(\frac{B(3, Y_n^+(t_i)+1)}{B(1, Y_n^+(t_i)+1)} - \left(\frac{B(2, Y_n^+(t_i)+1)}{B(1, Y_n^+(t_i)+1)}\right)^2\right). \end{aligned}$$

Lemma 2 together with (12) in Lemma 3 and (13) in Lemma 4 yields $\sup_{i=1,\dots,q_n} |W_{ni}| = O(n^{-2-\alpha})$. Hence

$$\mathrm{Var}(Z_n(t) - Z_n(s)|\mathrm{D}_n)$$

$$= \sum_{s < t_i \leq t} n[\mathrm{E}((\Delta A_d(t_i))^2 | \mathrm{D}_n) - (\mathrm{E}(\Delta A_d(t_i) | \mathrm{D}_n))^2]$$

$$= \int_s^t \frac{n}{Y_n(u)} \left[ \frac{Y_n^2(u)(Y_n^+(u)+1)}{(Y_n^+(u)+2)^2(Y_n^+(u)+3)} \right] \frac{dN_n(u)}{Y_n(u)} + \sum_{s < t_i \leq t} n W_{ni}.$$

Since

$$\sup_{u \in [0,\tau]} \left| \left[ \frac{Y_n^2(u)(Y_n^+(u)+1)}{(Y_n^+(u)+2)^2(Y_n^+(u)+3)} \right] - 1 \right| \to 0$$

and

$$\sup_{u \in [0,\tau]} \left| \frac{n}{Y_n(u)} - Q(u)^{-1} \right| \to 0$$

with probability 1, we have by Lemma 6,

$$\int_s^t \frac{n}{Y_n(u)} \left[ \frac{Y_n^2(u)(Y_n^+(u)+1)}{(Y_n^+(u)+2)^2(Y_n^+(u)+3)} \right] \frac{dN_n(u)}{Y_n(u)} \to U_0(t) - U_0(s)$$

uniformly in $s$ and $t$ with probability 1. Since

$$\sum_{s < t_i \leq t} n|W_{ni}| \leq n^2 \sup_{i=1,\ldots,q_n} |W_{ni}| = O(n^{-\alpha}),$$

we obtain

(19) $$\sup_{s,t \in [0,\tau]} |\mathrm{Var}(Z_n(t) - Z_n(s) | \mathrm{D}_n) - (U_0(t) - U_0(s))| \to 0$$

with probability 1. Now (18) and (19) imply the convergence of the finite dimensional posterior distributions of $Z_n$ to those of $W(U_0)$ with probability 1.

Finally, note that

$$\Pr\{|Z_n(t) - Z_n(s)| \geq \varepsilon | \mathrm{D}_n\} \leq \frac{1}{\varepsilon^2} \mathrm{Var}(Z_n(t) - Z_n(s) | \mathrm{D}_n).$$

By (19), we have

$$\mathrm{Var}(Z_n(t) - Z_n(s) | \mathrm{D}_n) = U_0(t) - U_0(s) + o(1)$$

with probability 1. Since $U_0(t)$ is continuous, with probability 1 we can make $\Pr\{|Z_n(t) - Z_n(s)| \geq \varepsilon | \mathrm{D}_n\}$ as small as possible for sufficiently large $n$ by choosing $t$ and $s$ sufficiently close. Hence by Theorem 19 in Section V.4 in Pollard (1984) we conclude that $Z_n$ given $\mathrm{D}_n$ converges weakly to $W(U_0)$ on $D[0,\tau]$ with probability 1.  □

PROOF OF THEOREM 1(iii). Let $W_{ni} = \mathrm{E}(\Delta A_d(t_i)|\mathrm{D}_n) - 1/Y_n^+(t_i)$. Then Lemma 4 yields $\sup_{i=1,\dots,q_n} |W_{ni}| = O(n^{-1-\alpha})$. Since

$$\mathrm{E}(A_d(t)|\mathrm{D}_n) = \int_0^t \frac{Y_n(s)}{Y_n^+(s)} \frac{dN_n(s)}{Y_n(s)} + \sum_{t_i \le t} W_{ni},$$

we have

$$(20) \qquad \sup_{t\in[0,\tau]} |\mathrm{E}(A_d(t)|\mathrm{D}_n) - \hat{A}_n(t)| \le \int_0^\tau \Big|1 - \frac{Y_n(s)}{Y_n^+(s)}\Big| \frac{dN_n(s)}{Y_n(s)} + O(n^{-\alpha}).$$

Since the first term on the right-hand side of (20) is $O(n^{-1})$ by Lemma 6, the proof is done. □

## APPENDIX A

**Proving Lemma 1.** Let

$$B_\alpha(s) = \int_0^1 x^\alpha (1-x)^{Y_n^+(s)}\, dx = \frac{\Gamma(\alpha+1)\Gamma(Y_n^+(s)+1)}{\Gamma(Y_n^+(s)+\alpha+2)}.$$

Then Lemma 2 yields that

$$(21) \qquad \sup_{s\in[0,\tau]} |(Y_n^+(s)+1)^{\alpha+1} B_\alpha(s)| \to \Gamma(\alpha+1)$$

and

$$(22) \qquad \sup_{s\in[0,\tau]} |(Y_n^+(s)+1) B_\alpha(s)| = O(n^{-\alpha})$$

with probability 1 for $\alpha > 0$.

LEMMA 5. *With probability 1, we have:*

(i)

$$(23) \qquad \sup_{s\in[0,\tau]} \Big|\frac{Y_n(s)B_1(s)}{B_0(s)+B_\alpha(s)} - 1\Big| = O(n^{-\min\{1,\alpha\}});$$

(ii) *if $\alpha \le 1/2$,*

$$(24) \qquad \sup_{s\in[0,\tau]} \Big|n^\alpha\Big(\frac{Y_n(s)B_1(s)}{B_0(s)+B_\alpha(s)} - 1\Big) + \Gamma(\alpha+1)Q(s)^{-\alpha}\Big| \to 0;$$

(iii)

$$(25) \qquad \sup_{s\in[0,\tau]} \Big|\frac{n^\alpha Y_n(s)B_{\alpha+1}(s)}{B_0(s)+B_\alpha(s)} - \Gamma(\alpha+2)Q(s)^{-\alpha}\Big| \to 0$$

*for $\alpha > 0$.*

PROOF. For (23), (22) yields

$$\left| \frac{Y_n(s)B_1(s)}{B_0(s) + B_\alpha(s)} - 1 \right|$$

$$= \left| \frac{Y_n(s) - Y_n^+(s) - 2}{(Y_n^+(s) + 2)(1 + (Y_n^+(s) + 1)B_\alpha(s))} - \frac{(Y_n^+(s) + 1)B_\alpha(s)}{1 + (Y_n^+(s) + 1)B_\alpha(s)} \right|$$

$$\leq \frac{2}{Y_n(\tau)} + \sup_{s \in [0,\tau]} |(Y_n^+(s) + 1)B_\alpha(s)|$$

$$= O(n^{-1}) + O(n^{-\alpha})$$

$$= O(n^{-\min\{1,\alpha\}})$$

with probability 1.

For (24), note that

$$n^\alpha \left( \frac{Y_n(s)B_1(s)}{B_0(s) + B_\alpha(s)} - 1 \right) + \Gamma(\alpha+1)Q(s)^{-\alpha}$$

$$(26) \qquad = \frac{n^\alpha(Y_n(s) - Y_n^+(s) - 2)}{(Y_n^+(s) + 2)(1 + (Y_n^+(s) + 1)B_\alpha(s))}$$

$$(27) \qquad - \left[ \frac{n^\alpha(Y_n^+(s) + 1)B_\alpha(s)}{1 + (Y_n^+(s) + 1)B_\alpha(s)} - \Gamma(\alpha+1)Q(s)^{-\alpha} \right].$$

Since $\alpha \leq 1/2$, $\sup_{s \in [0,\tau]} |(26)| \leq 2n^\alpha/Y_n(\tau) \to 0$ with probability 1. For (26), let $p(s) = \Gamma(\alpha+1)Q(s)^{-\alpha}$. Then

$$|(27)| \leq |n^\alpha(Y_n^+(s) + 1)B_\alpha(s) - p(s)| + |p(s)(Y_n^+(s) + 1)B_\alpha(s)|.$$

Here

$$|n^\alpha(Y_n^+(s) + 1)B_\alpha(s) - p(s)|$$

$$\leq \Gamma(\alpha+1) \frac{(Y_n^+(s) + 1)^{\alpha+1}\Gamma(Y_n^+(s) + 1)}{\Gamma(Y_n^+(s) + \alpha + 2)}$$

$$\times \left| \left( \frac{n}{Y_n^+(s) + 1} \right)^\alpha - Q(s)^{-\alpha} \right|$$

$$+ \Gamma(\alpha+1)Q(s)^{-\alpha} \left| \frac{(Y_n^+(s) + 1)^{\alpha+1}\Gamma(Y_n^+(s) + 1)}{\Gamma(Y_n^+(s) + \alpha + 2)} - 1 \right|.$$

Since $\sup_{t \in [0,\tau]} Y_n^+(t) = O(n)$, we conclude $\sup_{s \in [0,\tau]} |n^\alpha(Y_n^+(s) + 1)B_\alpha(s) - p(s)| \to 0$ with probability 1 by (9). Also we have $\sup_{s \in [0,\tau]} |p(s)(Y_n^+(s) + 1)B_\alpha(s)| \to 0$ with probability 1 by (22) and the proof is done.

For (25), (21) yields

$$\frac{n^\alpha Y_n(s) B_{\alpha+1}(s)}{B_0(s) + B_\alpha(s)}$$

$$= \frac{(n/Y_n(s))^\alpha}{1 + (Y_n^+(s) + 1) B_\alpha(s)} \left(\frac{Y_n(s)}{Y_n^+(s) + 1}\right)^{1+\alpha} (Y_n^+(s) + 1)^{2+\alpha} B_{\alpha+1}(s)$$

$$\to \Gamma(\alpha+2) Q(s)^{-\alpha}$$

uniformly in $s \in [0, \tau]$ with probability 1. $\square$

PROOF OF LEMMA 1. Note that

$$\mathrm{E}(A_d(t)|\mathrm{D}_n) = \int_0^t \frac{1}{B_0(s) + B_\alpha(s)} \int_0^1 x(1-x)^{Y_n^+(s)}(1+x^\alpha)\,dx\,dN_n(s).$$

Hence, we have

$$\mathrm{E}(A_d(t)|D_n) - \hat{A}_n(t)$$

$$= \int_0^t \frac{B_1(s) + B_{\alpha+1}(s)}{B_0(s) + B_\alpha(s)}\,dN_n(s) - \hat{A}_n(t)$$

$$= \int_0^t \left(\frac{Y_n(s) B_1(s)}{B_0(s) + B_\alpha(s)} - 1\right)\frac{dN_n(s)}{Y_n(s)} + \int_0^t \frac{Y_n(s) B_{\alpha+1}(s)}{B_0(s) + B_\alpha(s)}\frac{dN_n(s)}{Y_n(s)}.$$

For $0 < \alpha \le 1/2$, (24) in Lemma 5 and Lemma 6 yield

$$(28) \qquad \sup_{t \in [0,\tau]} \left| \int_0^t n^\alpha \left(\frac{Y_n(s) B_1(s)}{B_0(s) + B_\alpha(s)} - 1\right)\frac{dN_n(s)}{Y_n(s)} \right.$$
$$\left. + \int_0^t \Gamma(\alpha+1) Q(s)^{-\alpha}\,dA_0(s) \right| \to 0$$

with probability 1, and (25) in Lemma 5 and Lemma 6 imply

$$(29) \qquad \sup_{t \in [0,\tau]} \left| \int_0^t n^\alpha \frac{Y_n(s) B_{\alpha+1}(s)}{B_0(s) + B_\alpha(s)}\frac{dN_n(s)}{Y_n(s)} \right.$$
$$\left. - \int_0^t \Gamma(\alpha+2) Q(s)^{-\alpha}\,dA_0(s) \right| \to 0$$

with probability 1. Combining (28) and (29), we have

$$\sup_{t \in [0,\tau]} \left| n^\alpha (\mathrm{E}(A_d(t)|\mathrm{D}_n) \right.$$
$$\left. - \hat{A}_n(t)) - \int_0^t (\Gamma(\alpha+2) - \Gamma(\alpha+1)) Q(s)^{-\alpha}\,dA_0(s) \right| \to 0$$

with probability 1.

For $\alpha > 1/2$, (23) in Lemma 5 and Lemma 6 yield

$$(30) \qquad \sup_{t \in [0,\tau]} \left| \int_0^t n^{1/2} \left( \frac{Y_n(s)B_1(s)}{B_0(s) + B_\alpha(s)} - 1 \right) \frac{dN_n(s)}{Y_n(s)} \right| \to 0$$

with probability 1, and (25) in Lemma 5 and Lemma 6 imply

$$(31) \qquad \sup_{t \in [0,\tau]} \left| \int_0^t n^{1/2} \frac{Y_n(s)B_{\alpha+1}(s)}{B_0(s) + B_\alpha(s)} \frac{dN_n(s)}{Y_n(s)} \right| \to 0$$

with probability 1. Combining (30) and (31), we have

$$\sup_{t \in [0,\tau]} |n^{1/2}(\mathrm{E}(A_d(t)|\mathrm{D}_n) - \hat{A}_n(t))| \to 0$$

with probability 1.  $\square$

## APPENDIX B

**Technical lemmas.**

LEMMA 6. *Let $X_1(t), X_2(t), \ldots$ be stochastic processes defined on $[0, \tau]$. Suppose that there exists a continuous function $X(t)$ defined on $[0, \tau]$ such that*

$$\lim_{n \to \infty} \sup_{t \in [0,\tau]} |X_n(t) - X(t)| = 0$$

*with probability 1. Then*

$$\sup_{t \in [0,\tau]} \left| \int_0^t X_n(s) \frac{1}{Y_n(s)} \, dN_n(s) - \int_0^t X(s) \, dA_0(s) \right| \to 0$$

*with probability 1.*

PROOF. This lemma is an easy consequence of the Glivenko–Cantelli theorem and Lemma A.2 in Tsiatis (1981).  $\square$

LEMMA 7. *Let $X_n$ be a sequence of subordinators such that $\mathrm{E}(X_n(t)) \to X_0(t)$ and $\mathrm{Var}(X_n(t)) \to 0$ for some continuous function $X_0(t)$ and all $t \in [0, \tau]$. Then $\mathcal{L}(X_n) \xrightarrow{d} \delta_{X_0}$ on $D[0, \tau]$.*

PROOF. Note that $X_0$ should be a monotonically increasing function since $X_n$ are subordinators. Hence, the continuity of $X_0$ together with the assumptions implies that $\sup_{t \in [0,\tau]} |X_n(t) - X_0(t)| \to 0$ in probability.  $\square$

Department of Statistics
Seoul National University
Seoul 151-742
Korea
e-mail: ydkim0903@naver.com
e-mail: leej@stats.snu.ac.kr